\documentclass{amsart}
\usepackage{a4}		
\usepackage{amsmath}
\usepackage{amssymb}

\newtheorem{theorem}{Theorem}[section]

\newtheorem{proposition}[theorem]{Proposition}

\theoremstyle{definition}

\theoremstyle{remark}

\numberwithin{equation}{section}

\def\Q{{\mathbb{Q}}}\def\Z{{\mathbb{Z}}}
\def\norm{{\mathcal{N}}}\def\trace{{\mathcal{T}}}
\def\divides{\mid}

\begin{document}

\title{Elliptic curves with good reduction away from 2: III}
\author{Richard G.E. Pinch}
\address{Queens' College, University of Cambridge, Silver Street, Cambridge CB3 9ET, U.K.}
\email{rgep@cam.ac.uk}


\date{4 March 1998}

\begin{abstract}
We list the elliptic curves defined over $\Q ( \sqrt 5 )$ 
with good reduction away from 2.  There are 368 isomorphism classes.  
We give a global minimal model for each class. 
\end{abstract}

\maketitle

\section{Introduction }

In this paper we list the elliptic curves defined over $\Q ( \sqrt 5 )$ 
with good reduction away from 2.  We use the results of the previous papers
\cite{Pin:egr2} and \cite{Pin:egr2ii}, referred to as I and II
respectively.  By Theorem 1.14 of I, such a curve must have a point of order
2 defined over $\Q ( \sqrt 5 )$ and by Theorem 2.3 of II, 
if $t \in \Q ( \sqrt 5 )$ is the
corresponding value of the Hauptmodul on $X_0 (2)$ then either $t$ or
$t' = 4096 / t$
satisfies one of the equations
\begin{equation}\label{eq01}
t = 64 u / v, \qquad u + v = x^2
\end{equation}
or
\begin{equation}\label{eq02}
t = 64v / 2^a u, \qquad 2^a u + v = x^2
\end{equation}
where $u$, $v$ are units, $x \in \Q ( \sqrt 5 )$ and $a \ge 0$.
We solve these equations and determine the corresponding $j$-invariants
and obtain a global minimal equation in each isomorphism class by Tate's
algorithm.  There are 368 isomorphism classes of these curves.

\section{Diophantine equations over the quadratic field}

In this section we solve the Diophantine equations (\ref{eq01}) and (\ref{eq02}) 
over $\Q ( \sqrt 5 )$.  The fundamental unit 
is $\epsilon = \frac{1 + \sqrt 5}{2}$ and the ring of integers is
$\Z [ \epsilon ] $.  
Let $\norm$ and $\trace$ denote absolute norm and trace: let $'$ denote the
non-trivial automorphism of $\Q ( \sqrt 5 )$.

\begin{theorem}\label{thm11}
The equation
$$
2 x^2 = u + v,  \qquad u, v  \mbox{ units of } \Q ( \sqrt 5 ) ,
$$
has the solutions
$$
\begin{array}{rrrr}
\mbox{a)} &  x = 0 &  u = 1    &  v = -1           \\
\mbox{b)} & 1          & 1             & 1                     \\
\mbox{c)} & 1          & \epsilon^2 & \epsilon'         \\
\mbox{d)} & 3          & \epsilon^6 & {\epsilon'}^6    \\
\end{array}
$$
and all other solutions are obtained by conjugation and scaling by units.

The equations
$$
x^2 = 2^a u + v, \qquad u, v \mbox{ units,} \quad  a \ge 0,
$$
have the solutions
$$
\begin{array}{rrrrr}
\mbox{e)} & a=0 &     x =      0 & u=       1 & v= -1 \\
\mbox{f)} &   0 &              1 & \epsilon   & \epsilon' \\
\mbox{g)} &   0 &              1 & \epsilon^2 & - \epsilon \\
\mbox{h)} &   0 &              2 & \epsilon^3 & {\epsilon'}^3 \\
\mbox{i)} &   1 &              1 &          1 & -1 \\
\mbox{j)} &   1 &              1 & - \epsilon & \epsilon^3 \\
\mbox{k)} &   1 &       \epsilon &          1 & - \epsilon' \\
\mbox{l)} &   1 &        \sqrt 5 & \epsilon^2 & {\epsilon'}^3 \\
\mbox{m)} &   1 & 8 + 15 \epsilon & \epsilon^{13} & \epsilon' \\
\mbox{n)} &   2 &        \sqrt 5 &          1 & 1 \\
\mbox{o)} &   2 &     \epsilon^3 & \epsilon^3 & 1 \\
\mbox{p)} &   3 &              3 &          1 & 1 \\
\mbox{q)} &   3 & 5 + 2 \sqrt 5  & \epsilon^5 & 1 \\
\mbox{r)} &   3 & 3 + 2 \sqrt 5  & \epsilon^4 & 1 \\
\mbox{s)} &   4 & 17 + 8 \sqrt 5 & \epsilon^9 & 1  \\
\end{array}
$$
and all other solutions are obtained by conjugation or scaling by units.

The corresponding values of $t$ are
$$
\begin{array}{rr}
\mbox{a,e)} &    - 64   \\
\mbox{b)} &    64   \\
\mbox{c)} &    - 64 \epsilon^3 \\
\mbox{d,h)} &    - 64 \epsilon^6 \\
\mbox{f)} &    - 64 \epsilon^2 \\
\mbox{g)} &    64 \epsilon   \\
\mbox{i)} &    - 32   \\
\mbox{j)} &    - 32 \epsilon   \\
\mbox{k)} &    - 32 \epsilon^2 \\
\mbox{l)} &    32 \epsilon^5 \\
\mbox{m)} &    - 32 \epsilon^{14} \\
\mbox{n)} &    16   \\
\mbox{o)} &    - 16 \epsilon^3 \\
\mbox{p)} &    8    \\
\mbox{q)} &    - 8 \epsilon^5 \\
\mbox{r)} &    8 \epsilon^4 \\
\mbox{s)} &   - 4 \epsilon^9 \\
\end{array}
$$
and their conjugates.
\end{theorem}
\begin{proof}
The Theorem will follow from Propositions \ref{prop13}, \ref{prop15} 
and \ref{prop16}.
\end{proof}

We define  {\em  Lucas numbers} $L_n$ and  {\em  Fibonacci numbers} $F_n$ 
to satisfy the recurrence relation
$$
X_{{n+1}} = X_n + X_{{n-1}}
$$
with the initial conditions $L_0 = 2$, $L_1 = 1$ and
$F_0 = 0$, $F_1 = 1$ respectively.  Equivalently, we have
$$
\epsilon^n = \frac{L_n + F_n \sqrt 5}{2}
$$

We quote the next result from Cohn \cite{Coh:LucFib} Theorems 1 and 2.

\begin{proposition}\label{prop12}
The Diophantine equations
\begin{eqnarray*}
      L_n & = & X^2 ; \\
      L_n & = & 2 X^2 \\
\end{eqnarray*}
have the solutions
$$
n = 1,3,     \qquad X = \pm 1, \pm 3; \\
$$
and
$$
n=0, \pm  6, \qquad X = \pm 1, \pm 3  \\
$$
respectively.
\end{proposition}

\begin{proposition}\label{prop13}
The equations
\begin{eqnarray*}
     2 x^2 & = & u + v ; \\
      x^2  & = & u + v   \\
\end{eqnarray*}
with $u$, $v$ units, have the solutions
\begin{eqnarray*}
             2.0^2 = & 1 - 1           \\
               2.1^2 = & 1 + 1                    \\
               2.1^2 = & \epsilon^2 + \epsilon'   \\
               2.3^2 = & \epsilon^6 + {\epsilon'}^6 \\
\end{eqnarray*}
and
\begin{eqnarray*}
             0^2 =   & 1 - 1           \\
               1^2 = & \epsilon + \epsilon'     \\
               1^2 = & \epsilon^2 - \epsilon    \\
               2^2 = & \epsilon^3 + {\epsilon'}^3 \\
\end{eqnarray*}
respectively.
\end{proposition}

\begin{proof}
Let $A$ be 1 or 2, according to the equation we consider.
Put $U = \norm u$, $V = \norm v$, $X = \norm x$, $w = u v'$.
Taking norms we have
$$
A^2 X^2 = U + V + \trace w
$$

If $U=-V$ then $\trace w = (AX)^2$, so that by Proposition \ref{prop12} we have
$AX = \pm 1$ or $\pm 2$. But 2 is not a norm, so $X$ must be $1$, giving
$w=\epsilon$ or $\epsilon^3$, and solutions
$1^2 = \epsilon^2 - \epsilon$ and
$2.1^2 = \epsilon^2 + \epsilon'$.

If $U=V$, then $u^{{-1}} = U u'$, $v^{{-1}} = V v'$.
Now
$$
A {x'}^2 = u' + v' =
       U \left( u^{{-1}} + v^{{-1}} \right) = \frac{U}{uv} (u+v)
        = \frac{U}{uv} A x^2 ,
$$
so that
$$
{x'}^2 = \frac{U}{uv} x^2 ,
$$
and, considering the factorisation of $x^2$ into prime elements,
we have $x^2 = Y w$ where $Y \in \Z$ and $w$ is a unit.
Scaling by $w$, we may take $x = Z$ or $Z \sqrt 5$ where
$Z \in \Z$. But now $x^2 \in \Z$
and so $u - u' = v' - v$.  Thus either $u = \pm 1$ and
$v = \pm 1$, or $u = v'$.  The first case gives us the solutions
$A.0^2 = 1 - 1$ and $2.1^2 = 1 + 1$ and the second gives
the equation $A x^2 = \pm L_n$ for some $n$.  Considering signs,
$A x^2 = L_n$.  We cannot have $x^2 = 5 Z^2$, for
then $\sqrt 5$ would divide $\epsilon^n$; so $x^2 = Z^2$
and we have $A Z^2 = L_n$.  Applying Proposition \ref{prop12} we have the
remaining solutions.
\end{proof}

The next result is established in proposition 9 of \cite{Pin:sPell} using 
techniques of transcendence theory: the method is further developed in \cite{Pin:squad}.

\begin{proposition}\label{prop14}
The pairs of simultaneous Diophantine equations
\begin{eqnarray*}
     (X+2)^2 - 10Y^2 & = & -1 ,          \\
     (X-2)^2 - 2 Z^2 & = & -1 ;          \\
\end{eqnarray*}
and
\begin{eqnarray*}
      X^2 - 2 Y^2  & = & -1 ,            \\
      X^2 - 10 Z^2 & = & -9              \\
\end{eqnarray*}
have the solutions
$$
     (X,Y,Z) = (1, \pm 1, \pm 1)  \mbox{ or }  (-5, \pm 1, \pm 5) ; 
$$
and
$$
      (X,Y,Z) = ( \pm 1, \pm 1, \pm 1)  \mbox{ or }  
                          ( \pm 41, \pm 29, \pm 13) 
$$
respectively.
\end{proposition}

\begin{proposition}\label{prop15}
The equation
$$
x^2 = 2u + v,  \qquad u, v  \mbox{ units,}
$$
has the solutions
\begin{eqnarray*}
         1^2 & = &                  2.1 - 1                             \\
         1^2 & = &                  2(- \epsilon ) + \epsilon^3         \\
         \epsilon^2 & = &            2.1 - \epsilon'                    \\
         ( \sqrt 5 )^2 & = &         2 \epsilon^2 + {\epsilon'}^3         \\
         (8+15 \epsilon )^2 & = &   2 \epsilon^13 + \epsilon' .         \\
\end{eqnarray*}
\end{proposition}

\begin{proof}
Let $U = \norm u$, $V = \norm v$, $w = u v' = \frac{S + P \sqrt 5}{ 2 }$,
$W = \norm w$, $X = \norm x$.  We have $W = UV$,
$$
X^2 = 4 U + V + 2 S
$$
and
$$
S^2 - 5 P^2 = 4W .
$$
Clearly $X$ is odd.  If $X = 1$, then $x$ is a unit, and we solve the
equation by Proposition \ref{prop13}, giving the solutions
$1^2 = 2.1 - 1$, $1^2 = 2(- \epsilon ) + \epsilon^3$
and $\epsilon^2 = 2.1 - \epsilon'$.
We note that $X$ cannot be 3, which is not a norm from $\Q ( \sqrt 5 )$.

Assume now that $X \ge 5$, and
consider four cases according to the values of $U$ and $V$.

Suppose that $U = V = +1$.  Then $W = +1$ and, substituting, we have
$$
\left( \frac{X^2 - 5}{ 2 } \right)^2 - 5 P^2 = 4 ,
$$
so
$$
( X^2 - 1 )( X^2 - 9 ) = 5(2P)^2 .
$$
Since $X$ is odd, the two terms on the LHS are congruent to 0 and 8
$\mod 16$, so their highest common factor is 8; since $X \ge 5$,
each term is positive.   We must have either
\begin{eqnarray*}
X^2 - 1 & = & 10 Y^2 ,  \\
X^2 - 9 & = &  2 Z^2 ;  \\
\end{eqnarray*}
or
\begin{eqnarray*}
X^2 - 1 & = & 2 Y^2 ,  \\
X^2 - 9 & = & 10 Z^2 .  \\
\end{eqnarray*}
The first pair have no solution: since 2 is not a quadratic residue of 3,
3 must divide $X$, but then $-1$ would be a quadratic residue of 3, a
contradiction.  The second pair lead to no solution
of the original equation; for, from the second equation of the pair,
$X^2 \equiv 4 \mod (5)$, so $X \equiv \pm 2 \mod (5)$,
but this is impossible since $X$ is a norm from $\Q ( \sqrt 5 )$.  
Hence there is no solution in the case $U = V = +1$.

Suppose now that $U = +1$, $V = -1$, so that $W = -1$.
We have
$$
\left( \frac{X^2 - 3}{ 2 } \right)^2 - 5 P^2 = -4,
$$
so
$$
\left((X+2)^2 + 1\right)\left((X-2)^2 + 1\right) = 5(2P)^2 .
$$
Since $X$ is odd, each term on the LHS is congruent to $2 \mod (4)$ and so
their highest common factor is 2; further, each is positive.  So
\begin{eqnarray*}
(X \pm 2)^2 + 1 & = & 10 Y^2 ,  \\
(X \mp 2)^2 + 1 & = & 2 Z^2 .  \\
\end{eqnarray*}
By Proposition \ref{prop14}, the solutions are $(X,Y,Z) = (1, \pm 1, \pm 1)$,
dealt with above, and $(5, \pm 1, \pm 5)$, leading to
$( \sqrt 5 )^2 = 2 \epsilon^2 + {\epsilon'}^3$.

Suppose now that $U = -1$, $V = +1$, so that $W = -1$.  Substituting,
$$
\left( \frac{X^2 + 3}{ 2 } \right)^2 - 5 P^2 = -4
$$
so
$$
((X+1)^2 + 4)((X-1)^2 + 4) = 5(2P)^2 .
$$
Since $X$ is odd, each term on the LHS is divisible by 4.  The highest common
factor of the terms divides their difference, which is $4X$, and so is
either 4 (if $X$ is prime to 5) or 20 (if $5 \divides X$).  In either case, we
have
\begin{eqnarray*}
(X \pm 1)^2 + 4 & = & 5Y^2 ,  \\
(X \mp 1)^2 + 4 & = & Z^2 .  \\
\end{eqnarray*}
The second equation implies that $(X,Z) = (2, \pm 2)$ or $(0, \pm 2)$; 
but $X$ must be odd, so there is no solution in this case.

Finally, suppose that $U = V = -1$, so that $W = +1$.  Substituting,
$$
\left( \frac{X^2 + 5}{ 2 } \right)^2 - 5 P^2 = 4 ,
$$
so
$$
(X^2 + 1)(X^2 + 9) = 5 (2P)^2 .
$$
Since $X$ is odd, each term on the LHS is congruent to $2 \mod (4)$ and so
their highest common factor is 2; further, each is positive.
We must have either
\begin{eqnarray*}
X^2 + 1 & = & 10 Y^2 ,   \\
X^2 + 9 & = & 2 Z^2   \\
\end{eqnarray*}
or
\begin{eqnarray*}
X^2 + 1 & = & 2 Y^2 ,   \\
X^2 + 9 & = & 10 Z^2 .   \\
\end{eqnarray*}
In the first case, $X^2 \equiv 4 \mod (5)$, so that $X \equiv \pm 2 \mod (5)$:
but $X$ is a norm from $\Q ( \sqrt 5 )$ and this is impossible.  
In the second case,
by Proposition \ref{prop14}, the solutions are $(X,Y,Z) = ( \pm 1, \pm 1, \pm 1 )$,
dealt with above, and $( \pm 41, \pm 29, \pm 13)$ giving the solution
$x = 8+15 \epsilon$, $u = \epsilon^{13}$, $v = \epsilon'$.
\end{proof}

\begin{proposition}\label{prop16}
The equation
$$
x^2 = 2^a u + v , \qquad a \ge 2, \qquad u, v \mbox{ units,}
$$
has solutions
$$
\begin{array}{rrrr}
 a = 2 &  x = \sqrt 5  &    u = 1 &    v = 1 \\
 2 &  \epsilon^3  &    \epsilon^3 &    1 \\
 3 &  3 &    1 &    1 \\
 3 &  5 + 2 \sqrt 5 &    \epsilon^5 &    1 \\
 3 &  3 + 2 \sqrt 5 &    \epsilon^4 &    1 \\
 4 &  17 + 8 \sqrt 5 &    \epsilon^9 &    1 . \\
\end{array}
$$
\end{proposition}

\begin{proof}
We have $v \equiv x^2 \mod (4)$, so $v$ is a square unit; scaling, we may
assume that $v = 1$.  Now
$$
x^2 - 1 = 2^a u ,
$$
and, factorising and choosing the sign on $x$, we have
\begin{eqnarray*}
x + 1 & = & 2^{{a-1}} u_1 ,  \\
x - 1 & = & 2 u_2 ,  \\
\end{eqnarray*}
where $u_1$, $u_2$ are units with $u = u_1 u_2$.
Subtracting,
$$
2^{{a-2}} u_1 - u_2 = 1,
$$
that is,
$$
2^{{a-2}} = u_1^{{-1}} + u_2 u_1^{{-1}} .
$$
This is solved by Proposition \ref{prop13}, taking the solutions with $x$ a power of 2,
and the result follows.
\end{proof}

\section{Equations of curves with good reduction away from 2 }

By Corollary 1.14.1 of II, an elliptic curve $E$ defined 
over $\Q ( \sqrt 5 )$ with good reduction
away from 2 has a point $P$ of order 2 defined over $\Q ( \sqrt 5 )$.  
By Proposition 2.1 of I, the corresponding invariant $j$ is given by
$$
j(t) = \frac{(t+16)^3}{ t }
$$
where $t \in \Q ( \sqrt 5 )$ is the value of the Hauptmodul for $X_0(2)$
corresponding to the pair $(E,P)$, and
the isogenous curve $E'$ has invariant $j' = j(t' )$, with
$t t' = 4096$.

The possible values of $t$ are obtained from Theorem \ref{thm11}.
Table 1 lists values of $t$ corresponding to a curve with good reduction away
from 2, together with the corresponding values of $j$ 
and their factorisations.
As in Section 5 of I we can list the minimal equation of a representative 
of each isomorphism class over $\Q ( \sqrt 5 )$, and these are given in
Table 2.  We express the curve in the form
$$
Y^2 + a_1 XY + a_3 Y =
                X^3 + a_2 X^2 + a_4 X  + a_6
$$
and give the coefficients $a_i$, the discriminant $\Delta$, the exponent of
the conductor $2^f$ and the reduction type according to Tate's algorithm
(see Tate \cite{Tate:letter}, Cremona \cite{Cre:book}, or \cite{Pin:egr}).
In each table {\em code} refers to the $j$ invariant.

\begin{theorem}
There are 56 values of $t$ and 43 values of $j$ corresponding to elliptic
curves over $\Q ( \sqrt 5 )$ with good reduction away from 2.
There are 368 isomorphism classes of such curves.
\end{theorem}

It is interesting to observe that every curve in Table 2 has additive bad reduction
In addition, we see that there is no curve with good reduction everywhere 
over $\Q ( \sqrt 5 )$: see Cremona \cite{Cre:Gamma1}, 
Ishii \cite{Ish:nonegr},\cite{Ish:nonegr2} and Pinch \cite{Pin:egr}.

\newpage

{\small
\def\e{\epsilon}
$$
\begin{array}{rlcrl}
  t
&  
&  code
&  j
&  
   \\
  -220+136\e
&  +\e^{-9}2^{2}
&   34 
&  78608
&  +2^{4}(17)^{3}
   \\
  -84-136\e
&  -\e^{9}2^{2}
&   34 
&  78608
&  +2^{4}(17)^{3}
   \\
  -64+40\e
&  +\e^{-5}2^{3}
&   28 
&  4928+960\e
&  +\e^{2}2^{6}(5-\e)^{3}
   \\
  40-24\e
&  +\e^{-4}2^{3}
&   31 
&  5888-960\e
&  +\e^{-2}2^{6}(4+\e)^{3}
   \\
  8
&  +2^{3}
&   11 
&  1728
&  +2^{6}(3)^{3}
   \\
  16+24\e
&  +\e^{4}2^{3}
&   28 
&  4928+960\e
&  +\e^{2}2^{6}(5-\e)^{3}
   \\
  -24-40\e
&  -\e^{5}2^{3}
&   31 
&  5888-960\e
&  +\e^{-2}2^{6}(4+\e)^{3}
   \\
  -48+32\e
&  +\e^{-3}2^{4}
&   14 
&  2048
&  +2^{11}
   \\
  16
&  +2^{4}
&   14 
&  2048
&  +2^{11}
   \\
  -16-32\e
&  -\e^{3}2^{4}
&   14 
&  2048
&  +2^{11}
   \\
  -19520+12064\e
&  -\e^{-14}2^{5}
&   40 
&  525604480-324907648\e
&  -\e^{-7}2^{7}(50-13\e)^{3}
   \\
  256-160\e
&  -\e^{-5}2^{5}
&   2 
&  103808-64640\e
&  -\e^{-7}2^{7}(4+\e)^{3}
   \\
  -64+32\e
&  -\e^{-2}2^{5}
&   17 
&  2688-1664\e
&  -\e^{-7}2^{7}
   \\
  -32+32\e
&  +\e^{-1}2^{5}
&   6 
&  1280+640\e
&  -\e2^{7}(1-2\e)^{3}
   \\
  -32
&  -2^{5}
&   5 
&  128
&  +2^{7}
   \\
  -32\e
&  -\e2^{5}
&   13 
&  1920-640\e
&  -\e^{-1}2^{7}(1-2\e)^{3}
   \\
  -32-32\e
&  -\e^{2}2^{5}
&   1 
&  1024+1664\e
&  +\e^{7}2^{7}
   \\
  96+160\e
&  +\e^{5}2^{5}
&   23 
&  39168+64640\e
&  +\e^{7}2^{7}(5-\e)^{3}
   \\
  -7456-12064\e
&  -\e^{14}2^{5}
&   41 
&  200696832+324907648\e
&  +\e^{7}2^{7}(37+13\e)^{3}
   \\
  -832+512\e
&  -\e^{-6}2^{6}
&   39 
&  914880-565760\e
&  +\e^{-12}2^{6}(-1-8\e)^{3}
   \\
  192-128\e
&  -\e^{-3}2^{6}
&   32 
&  63168-39040\e
&  -\e^{-15}2^{6}
   \\
  -192+128\e
&  +\e^{-3}2^{6}
&   26 
&  44864-26496\e
&  +\e^{-6}2^{6}(5+2\e)^{3}
   \\
  -128+64\e
&  -\e^{-2}2^{6}
&   8 
&  15040-9280\e
&  -\e^{-7}2^{6}(1-2\e)^{3}
   \\
  64-64\e
&  -\e^{-1}2^{6}
&   4 
&  12032-7232\e
&  +\e^{-8}2^{6}(3+2\e)^{3}
   \\
  -64+64\e
&  +\e^{-1}2^{6}
&   31 
&  5888-960\e
&  +\e^{-2}2^{6}(4+\e)^{3}
   \\
  -64
&  -2^{6}
&   11 
&  1728
&  +2^{6}(3)^{3}
   \\
  64
&  +2^{6}
&   35 
&  8000
&  +2^{6}(5)^{3}
   \\
  -64\e
&  -\e2^{6}
&   28 
&  4928+960\e
&  +\e^{2}2^{6}(5-\e)^{3}
   \\
  64\e
&  +\e2^{6}
&   27 
&  4800+7232\e
&  +\e^{8}2^{6}(5-2\e)^{3}
   \\
  -64-64\e
&  -\e^{2}2^{6}
&   30 
&  5760+9280\e
&  -\e^{7}2^{6}(1-2\e)^{3}
   \\
  -64-128\e
&  -\e^{3}2^{6}
&   12 
&  18368+26496\e
&  +\e^{6}2^{6}(7-2\e)^{3}
   \\
  64+128\e
&  +\e^{3}2^{6}
&   16 
&  24128+39040\e
&  +\e^{15}2^{6}
   \\
  -320-512\e
&  -\e^{6}2^{6}
&   22 
&  349120+565760\e
&  -\e^{12}2^{6}(9-8\e)^{3}
   \\
  -78080+48256\e
&  -\e^{-14}2^{7}
&   42 
&  8421373408-5204711200\e
&  +\e^{-13}2^{5}(-85+6\e)^{3}
   \\
  -1024+640\e
&  +\e^{-5}2^{7}
&   7 
&  1409888-870240\e
&  +\e^{-4}2^{5}(17-6\e)^{3}
   \\
  -256+128\e
&  -\e^{-2}2^{7}
&   33 
&  70368-43040\e
&  +\e^{-4}2^{5}(7-\e)^{3}
   \\
  128-128\e
&  -\e^{-1}2^{7}
&   24 
&  39680-22560\e
&  -\e2^{5}(9-8\e)^{3}
   \\
  -128
&  -2^{7}
&   3 
&  10976
&  +2^{5}(7)^{3}
   \\
  128\e
&  +\e2^{7}
&   10 
&  17120+22560\e
&  -\e^{-1}2^{5}(-1-8\e)^{3}
   \\
  -128-128\e
&  -\e^{2}2^{7}
&   18 
&  27328+43040\e
&  +\e^{4}2^{5}(6+\e)^{3}
   \\
  -384-640\e
&  -\e^{5}2^{7}
&   29 
&  539648+870240\e
&  +\e^{4}2^{5}(11+6\e)^{3}
   \\
  -29824-48256\e
&  -\e^{14}2^{7}
&   43 
&  3216662208+5204711200\e
&  -\e^{13}2^{5}(-79-6\e)^{3}
   \\
  768-512\e
&  -\e^{-3}2^{8}
&   36 
&  889584-548896\e
&  +\e^{-6}2^{4}(15-2\e)^{3}
   \\
  256
&  +2^{8}
&   34 
&  78608
&  +2^{4}(17)^{3}
   \\
  256+512\e
&  +\e^{3}2^{8}
&   20 
&  340688+548896\e
&  +\e^{6}2^{4}(13+2\e)^{3}
   \\
  4096-2560\e
&  -\e^{-5}2^{9}
&   15 
&  23528168-14540840\e
&  -\e^{-10}2^{3}(-29+5\e)^{3}
   \\
  2560-1536\e
&  +\e^{-4}2^{9}
&   38 
&  9036560-5578728\e
&  -\e^{-8}2^{3}(-34-3\e)^{3}
   \\
  512
&  +2^{9}
&   19 
&  287496
&  +2^{3}(33)^{3}
   \\
  1024+1536\e
&  +\e^{4}2^{9}
&   21 
&  3457832+5578728\e
&  -\e^{8}2^{3}(-37+3\e)^{3}
   \\
  1536+2560\e
&  +\e^{5}2^{9}
&   37 
&  8987328+14540840\e
&  -\e^{10}2^{3}(-24-5\e)^{3}
   \\
  56320-34816\e
&  -\e^{-9}2^{10}
&   25 
&  4386800300-2711191688\e
&  -\e^{-21}2^{2}(34-9\e)^{3}
   \\
  21504+34816\e
&  +\e^{9}2^{10}
&   9 
&  1675608612+2711191688\e
&  +\e^{21}2^{2}(25+9\e)^{3}
   \\
\end{array}
$$
}
Table 1. Invariants of elliptic curves defined over ${\mathbb Q}\left(\sqrt 5\right)$
with good reduction away from 2.

\newpage

{\small
\def\e{\epsilon}
$$
\begin{array}{cccccclrc}
 a_1 & a_2 & a_3 & a_4 & a_6 & \Delta & \mbox{type} & f & \mbox{code} \\
 0 & -\e & 0 & -16+11\e & 27-17\e & -\e^{-3}2^8 & I^*_1 & 3 & 36 \\
 0 & -\e & 0 & -171+106\e & -1050+647\e & -\e^32^{10} & III^* & 3 & 25 \\
 0 & -1+\e & 0 & -\e & 0 & +\e^62^4 & III & 3 & 14 \\
 0 & -1+\e & 0 & -5-6\e & -7-11\e & +\e^62^8 & I^*_1 & 3 & 34 \\
 0 & -1+\e & 0 & -5-11\e & 10+17\e & +\e^32^8 & I^*_1 & 3 & 20 \\
 0 & -1+\e & 0 & -65-106\e & -403-647\e & +\e^{-3}2^{10} & III^* & 3 & 9 \\
 0 & \e & 0 & -16+11\e & -27+17\e & -\e^{-3}2^8 & I^*_{0} & 4 & 36 \\
 0 & \e & 0 & -171+106\e & 1050-647\e & -\e^32^{10} & I^*_2 & 4 & 25 \\
 0 & 1 & 0 & -171-277\e & -1697-2747\e & +\e^32^{10} & I^*_2 & 4 & 9 \\
 0 & 1 & 0 & -448+277\e & -4444+2747\e & -\e^{-3}2^{10} & I^*_2 & 4 & 25 \\
 0 & 1+\e & 0 & \e & 0 & +2^4 & II & 4 & 14 \\
 0 & 1+\e & 0 & -4\e & -4-8\e & +\e^{-3}2^8 & I^*_{0} & 4 & 20 \\
 0 & 1+\e & 0 & -5+\e & -4+3\e & +2^8 & I^*_{0} & 4 & 34 \\
 0 & 1+\e & 0 & -5+6\e & 7-3\e & -\e^32^8 & I^*_{0} & 4 & 36 \\
 0 & 1-\e & 0 & -\e & 0 & +\e^62^4 & II & 4 & 14 \\
 0 & 1-\e & 0 & -5-6\e & 7+11\e & +\e^62^8 & I^*_{0} & 4 & 34 \\
 0 & 1-\e & 0 & -5-11\e & -10-17\e & +\e^32^8 & I^*_{0} & 4 & 20 \\
 0 & 1-\e & 0 & -65-106\e & 403+647\e & +\e^{-3}2^{10} & I^*_2 & 4 & 9 \\
 0 & -1 & 0 & -171-277\e & 1697+2747\e & +\e^32^{10} & I^*_2 & 4 & 9 \\
 0 & -1 & 0 & -448+277\e & 4444-2747\e & -\e^{-3}2^{10} & I^*_2 & 4 & 25 \\
 0 & -1-\e & 0 & \e & 0 & +2^4 & II & 4 & 14 \\
 0 & -1-\e & 0 & -4\e & 4+8\e & +\e^{-3}2^8 & I^*_{0} & 4 & 20 \\
 0 & -1-\e & 0 & -5+\e & 4-3\e & +2^8 & I^*_{0} & 4 & 34 \\
 0 & -1-\e & 0 & -5+6\e & -7+3\e & -\e^32^8 & I^*_{0} & 4 & 36 \\
 0 & 0 & 0 & 4 & 0 & -2^{12} & I^*_3 & 5 & 11 \\
 0 & 0 & 0 & 4+4\e & 0 & -\e^62^{12} & I^*_3 & 5 & 11 \\
 0 & 0 & 0 & -1 & 0 & +2^6 & III & 5 & 11 \\
 0 & 0 & 0 & -1-\e & 0 & +\e^62^6 & III & 5 & 11 \\
 0 & 0 & 0 & -11 & 14 & +2^9 & I^*_{0} & 5 & 19 \\
 0 & 0 & 0 & -11 & -14 & +2^9 & I^*_{0} & 5 & 19 \\
 0 & 0 & 0 & -11-11\e & 14+28\e & +\e^62^9 & I^*_{0} & 5 & 19 \\
 0 & 0 & 0 & -11-11\e & -14-28\e & +\e^62^9 & I^*_{0} & 5 & 19 \\
 0 & \e & 0 & -1 & -\e & +\e^22^6 & III & 5 & 31 \\
 0 & \e & 0 & -1-5\e & 3+2\e & +\e2^{12} & I^*_3 & 5 & 27 \\
 0 & \e & 0 & -11-10\e & -22-31\e & +\e^{-2}2^9 & I^*_{0} & 5 & 21 \\
 0 & \e & 0 & -32-53\e & -155-251\e & +\e^52^9 & I^*_{0} & 5 & 37 \\
 0 & -\e & 0 & -1 & \e & +\e^22^6 & III & 5 & 31 \\
 0 & -\e & 0 & -1-5\e & -3-2\e & +\e2^{12} & I^*_3 & 5 & 27 \\
 0 & -\e & 0 & -11-10\e & 22+31\e & +\e^{-2}2^9 & I^*_{0} & 5 & 21 \\
 0 & -\e & 0 & -32-53\e & 155+251\e & +\e^52^9 & I^*_{0} & 5 & 37 \\
 0 & 1 & 0 & 3-4\e & -5+4\e & +\e^{-5}2^{12} & I^*_3 & 5 & 27 \\
 0 & 1 & 0 & -1+4\e & -1-4\e & -\e^52^{12} & I^*_3 & 5 & 4 \\
 0 & 1 & 0 & -1-\e & -1-\e & +\e^42^6 & III & 5 & 28 \\
 0 & 1 & 0 & -2+\e & -2+\e & +\e^{-4}2^6 & III & 5 & 31 \\
 0 & 1 & 0 & -11-21\e & -37-59\e & +\e^{-1}2^9 & I^*_{0} & 5 & 37 \\
 0 & 1 & 0 & -32+21\e & -96+59\e & -\e2^9 & I^*_{0} & 5 & 15 \\
 0 & 1+\e & 0 & -21-31\e & -84-137\e & +\e^42^9 & I^*_{0} & 5 & 21 \\
 0 & 1+\e & 0 & -53+33\e & 168-105\e & +\e^{-4}2^9 & I^*_{0} & 5 & 38 \\
\end{array}
$$

$$
\begin{array}{cccccclrc}
 a_1 & a_2 & a_3 & a_4 & a_6 & \Delta & \mbox{type} & f & \mbox{code} \\
 0 & 1-\e & 0 & -1 & -1+\e & +\e^{-2}2^6 & III & 5 & 28 \\
 0 & 1-\e & 0 & -6+5\e & 5-2\e & -\e^{-1}2^{12} & I^*_3 & 5 & 4 \\
 0 & 1-\e & 0 & -21+10\e & -53+31\e & +\e^22^9 & I^*_{0} & 5 & 38 \\
 0 & 1-\e & 0 & -85+53\e & -406+251\e & -\e^{-5}2^9 & I^*_{0} & 5 & 15 \\
 0 & -1 & 0 & 3-4\e & 5-4\e & +\e^{-5}2^{12} & I^*_3 & 5 & 27 \\
 0 & -1 & 0 & -1+4\e & 1+4\e & -\e^52^{12} & I^*_3 & 5 & 4 \\
 0 & -1 & 0 & -1-\e & 1+\e & +\e^42^6 & III & 5 & 28 \\
 0 & -1 & 0 & -2+\e & 2-\e & +\e^{-4}2^6 & III & 5 & 31 \\
 0 & -1 & 0 & -11-21\e & 37+59\e & +\e^{-1}2^9 & I^*_{0} & 5 & 37 \\
 0 & -1 & 0 & -32+21\e & 96-59\e & -\e2^9 & I^*_{0} & 5 & 15 \\
 0 & -1+\e & 0 & -1 & 1-\e & +\e^{-2}2^6 & III & 5 & 28 \\
 0 & -1+\e & 0 & -6+5\e & -5+2\e & -\e^{-1}2^{12} & I^*_3 & 5 & 4 \\
 0 & -1+\e & 0 & -21+10\e & 53-31\e & +\e^22^9 & I^*_{0} & 5 & 38 \\
 0 & -1+\e & 0 & -85+53\e & 406-251\e & -\e^{-5}2^9 & I^*_{0} & 5 & 15 \\
 0 & -1-\e & 0 & -21-31\e & 84+137\e & +\e^42^9 & I^*_{0} & 5 & 21 \\
 0 & -1-\e & 0 & -53+33\e & -168+105\e & +\e^{-4}2^9 & I^*_{0} & 5 & 38 \\
 0 & 0 & 0 & \e & 0 & -\e^32^6 & II & 6 & 11 \\
 0 & 0 & 0 & 4\e & 0 & -\e^32^{12} & I^*_2 & 6 & 11 \\
 0 & 0 & 0 & -\e & 0 & +\e^32^6 & II & 6 & 11 \\
 0 & 0 & 0 & -4\e & 0 & +\e^32^{12} & I^*_2 & 6 & 11 \\
 0 & 0 & 0 & 1 & 0 & -2^6 & II & 6 & 11 \\
 0 & 0 & 0 & 1+\e & 0 & -\e^62^6 & II & 6 & 11 \\
 0 & 0 & 0 & 1-\e & 0 & +\e^{-3}2^6 & II & 6 & 11 \\
 0 & 0 & 0 & 4-4\e & 0 & +\e^{-3}2^{12} & I^*_2 & 6 & 11 \\
 0 & 0 & 0 & -1+\e & 0 & -\e^{-3}2^6 & II & 6 & 11 \\
 0 & 0 & 0 & -4 & 0 & +2^{12} & I^*_2 & 6 & 11 \\
 0 & 0 & 0 & -4+4\e & 0 & -\e^{-3}2^{12} & I^*_2 & 6 & 11 \\
 0 & 0 & 0 & -4-4\e & 0 & +\e^62^{12} & I^*_2 & 6 & 11 \\
 0 & 0 & 0 & -44 & 112 & +2^{15} & I^*_5 & 6 & 19 \\
 0 & 0 & 0 & -44 & -112 & +2^{15} & I^*_5 & 6 & 19 \\
 0 & 0 & 0 & -44-44\e & 112+224\e & +\e^62^{15} & I^*_5 & 6 & 19 \\
 0 & 0 & 0 & -44-44\e & -112-224\e & +\e^62^{15} & I^*_5 & 6 & 19 \\
 0 & \e & 0 & -\e & -1-\e & +\e2^6 & II & 6 & 27 \\
 0 & \e & 0 & -1-2\e & -2-3\e & +\e^52^6 & II & 6 & 30 \\
 0 & \e & 0 & -4-6\e & 4+6\e & +\e^32^6 & II & 6 & 22 \\
 0 & \e & 0 & -5-\e & -1+2\e & +\e^22^{12} & I^*_2 & 6 & 31 \\
 0 & \e & 0 & -5-9\e & 7+10\e & +\e^52^{12} & I^*_2 & 6 & 30 \\
 0 & \e & 0 & -17-25\e & -57-90\e & +\e^32^{12} & I^*_2 & 6 & 22 \\
 0 & \e & 0 & -45-41\e & 135+162\e & +\e^{-2}2^{15} & I^*_5 & 6 & 21 \\
 0 & \e & 0 & -65+43\e & 259-158\e & -\e^{-3}2^{14} & I^*_4 & 6 & 36 \\
 0 & \e & 0 & -129-213\e & 1027+1666\e & +\e^52^{15} & I^*_5 & 6 & 37 \\
 0 & \e & 0 & -685+423\e & -7977+4914\e & -\e^32^{16} & I^*_6 & 6 & 25 \\
 0 & -\e & 0 & -\e & 1+\e & +\e2^6 & II & 6 & 27 \\
 0 & -\e & 0 & -1-2\e & 2+3\e & +\e^52^6 & II & 6 & 30 \\
 0 & -\e & 0 & -4-6\e & -4-6\e & +\e^32^6 & II & 6 & 22 \\
 0 & -\e & 0 & -5-\e & 1-2\e & +\e^22^{12} & I^*_2 & 6 & 31 \\
 0 & -\e & 0 & -5-9\e & -7-10\e & +\e^52^{12} & I^*_2 & 6 & 30 \\
 0 & -\e & 0 & -17-25\e & 57+90\e & +\e^32^{12} & I^*_2 & 6 & 22 \\
\end{array}
$$

$$
\begin{array}{cccccclrc}
 a_1 & a_2 & a_3 & a_4 & a_6 & \Delta & \mbox{type} & f & \mbox{code} \\
 0 & -\e & 0 & -45-41\e & -135-162\e & +\e^{-2}2^{15} & I^*_5 & 6 & 21 \\
 0 & -\e & 0 & -65+43\e & -259+158\e & -\e^{-3}2^{14} & I^*_4 & 6 & 36 \\
 0 & -\e & 0 & -129-213\e & -1027-1666\e & +\e^52^{15} & I^*_5 & 6 & 37 \\
 0 & -\e & 0 & -685+423\e & 7977-4914\e & -\e^32^{16} & I^*_6 & 6 & 25 \\
 0 & 1 & 0 & \e & \e & -\e^52^6 & II & 6 & 4 \\
 0 & 1 & 0 & -\e & -\e & +\e^{-1}2^6 & II & 6 & 30 \\
 0 & 1 & 0 & 1-\e & 1-\e & +\e^{-5}2^6 & II & 6 & 27 \\
 0 & 1 & 0 & -1+\e & -1+\e & -\e2^6 & II & 6 & 8 \\
 0 & 1 & 0 & -1-4\e & -1+4\e & +\e^{-1}2^{12} & I^*_2 & 6 & 30 \\
 0 & 1 & 0 & -2-2\e & 2\e & +\e^{-3}2^6 & II & 6 & 22 \\
 0 & 1 & 0 & -4+2\e & 2-2\e & -\e^32^6 & II & 6 & 39 \\
 0 & 1 & 0 & -5+4\e & 3-4\e & -\e2^{12} & I^*_2 & 6 & 8 \\
 0 & 1 & 0 & -5-4\e & 3+4\e & +\e^42^{12} & I^*_2 & 6 & 28 \\
 0 & 1 & 0 & -9+4\e & 7-4\e & +\e^{-4}2^{12} & I^*_2 & 6 & 31 \\
 0 & 1 & 0 & -9-8\e & -9-24\e & +\e^{-3}2^{12} & I^*_2 & 6 & 22 \\
 0 & 1 & 0 & -17+8\e & -33+24\e & -\e^32^{12} & I^*_2 & 6 & 39 \\
 0 & 1 & 0 & -45-84\e & 251+388\e & +\e^{-1}2^{15} & I^*_5 & 6 & 37 \\
 0 & 1 & 0 & -129+84\e & 639-388\e & -\e2^{15} & I^*_5 & 6 & 15 \\
 0 & 1 & 0 & -685-1108\e & 12891+20868\e & +\e^32^{16} & I^*_6 & 6 & 9 \\
 0 & 1 & 0 & -1793+1108\e & 33759-20868\e & -\e^{-3}2^{16} & I^*_6 & 6 & 25 \\
 0 & 1+\e & 0 & -2+\e & -1 & +2^{10} & I^*_{0} & 6 & 14 \\
 0 & 1+\e & 0 & -2-19\e & 11+24\e & +\e^{-3}2^{14} & I^*_4 & 6 & 20 \\
 0 & 1+\e & 0 & -22+\e & 11-44\e & +2^{14} & I^*_4 & 6 & 34 \\
 0 & 1+\e & 0 & -22+21\e & -57+44\e & -\e^32^{14} & I^*_4 & 6 & 36 \\
 0 & 1+\e & 0 & -86-127\e & 459+756\e & +\e^42^{15} & I^*_5 & 6 & 21 \\
 0 & 1+\e & 0 & -214+129\e & -1429+884\e & +\e^{-4}2^{15} & I^*_5 & 6 & 38 \\
 0 & 1-\e & 0 & -1+\e & -2+\e & -\e^{-1}2^6 & II & 6 & 4 \\
 0 & 1-\e & 0 & -2-3\e & 1+2\e & +\e^62^{10} & I^*_{0} & 6 & 14 \\
 0 & 1-\e & 0 & -3+2\e & -5+3\e & -\e^{-5}2^6 & II & 6 & 8 \\
 0 & 1-\e & 0 & -6+\e & 1-2\e & +\e^{-2}2^{12} & I^*_2 & 6 & 28 \\
 0 & 1-\e & 0 & -10+6\e & 10-6\e & -\e^{-3}2^6 & II & 6 & 39 \\
 0 & 1-\e & 0 & -14+9\e & 17-10\e & -\e^{-5}2^{12} & I^*_2 & 6 & 8 \\
 0 & 1-\e & 0 & -22-23\e & -55-66\e & +\e^62^{14} & I^*_4 & 6 & 34 \\
 0 & 1-\e & 0 & -22-43\e & 101+158\e & +\e^32^{14} & I^*_4 & 6 & 20 \\
 0 & 1-\e & 0 & -42+25\e & -147+90\e & -\e^{-3}2^{12} & I^*_2 & 6 & 39 \\
 0 & 1-\e & 0 & -86+41\e & 297-162\e & +\e^22^{15} & I^*_5 & 6 & 38 \\
 0 & 1-\e & 0 & -262-423\e & -3063-4914\e & +\e^{-3}2^{16} & I^*_6 & 6 & 9 \\
 0 & 1-\e & 0 & -342+213\e & 2693-1666\e & -\e^{-5}2^{15} & I^*_5 & 6 & 15 \\
 0 & -1 & 0 & \e & -\e & -\e^52^6 & II & 6 & 4 \\
 0 & -1 & 0 & -\e & \e & +\e^{-1}2^6 & II & 6 & 30 \\
 0 & -1 & 0 & 1-\e & -1+\e & +\e^{-5}2^6 & II & 6 & 27 \\
 0 & -1 & 0 & -1+\e & 1-\e & -\e2^6 & II & 6 & 8 \\
 0 & -1 & 0 & -1-4\e & 1-4\e & +\e^{-1}2^{12} & I^*_2 & 6 & 30 \\
 0 & -1 & 0 & -2-2\e & -2\e & +\e^{-3}2^6 & II & 6 & 22 \\
 0 & -1 & 0 & -4+2\e & -2+2\e & -\e^32^6 & II & 6 & 39 \\
 0 & -1 & 0 & -5+4\e & -3+4\e & -\e2^{12} & I^*_2 & 6 & 8 \\
 0 & -1 & 0 & -5-4\e & -3-4\e & +\e^42^{12} & I^*_2 & 6 & 28 \\
 0 & -1 & 0 & -9+4\e & -7+4\e & +\e^{-4}2^{12} & I^*_2 & 6 & 31 \\
\end{array}
$$

$$
\begin{array}{cccccclrc}
 a_1 & a_2 & a_3 & a_4 & a_6 & \Delta & \mbox{type} & f & \mbox{code} \\
 0 & -1 & 0 & -9-8\e & 9+24\e & +\e^{-3}2^{12} & I^*_2 & 6 & 22 \\
 0 & -1 & 0 & -17+8\e & 33-24\e & -\e^32^{12} & I^*_2 & 6 & 39 \\
 0 & -1 & 0 & -45-84\e & -251-388\e & +\e^{-1}2^{15} & I^*_5 & 6 & 37 \\
 0 & -1 & 0 & -129+84\e & -639+388\e & -\e2^{15} & I^*_5 & 6 & 15 \\
 0 & -1 & 0 & -685-1108\e & -12891-20868\e & +\e^32^{16} & I^*_6 & 6 & 9 \\
 0 & -1 & 0 & -1793+1108\e & -33759+20868\e & -\e^{-3}2^{16} & I^*_6 & 6 & 25 \\
 0 & -1+\e & 0 & -1+\e & 2-\e & -\e^{-1}2^6 & II & 6 & 4 \\
 0 & -1+\e & 0 & -2-3\e & -1-2\e & +\e^62^{10} & I^*_{0} & 6 & 14 \\
 0 & -1+\e & 0 & -3+2\e & 5-3\e & -\e^{-5}2^6 & II & 6 & 8 \\
 0 & -1+\e & 0 & -6+\e & -1+2\e & +\e^{-2}2^{12} & I^*_2 & 6 & 28 \\
 0 & -1+\e & 0 & -10+6\e & -10+6\e & -\e^{-3}2^6 & II & 6 & 39 \\
 0 & -1+\e & 0 & -14+9\e & -17+10\e & -\e^{-5}2^{12} & I^*_2 & 6 & 8 \\
 0 & -1+\e & 0 & -22-23\e & 55+66\e & +\e^62^{14} & I^*_4 & 6 & 34 \\
 0 & -1+\e & 0 & -22-43\e & -101-158\e & +\e^32^{14} & I^*_4 & 6 & 20 \\
 0 & -1+\e & 0 & -42+25\e & 147-90\e & -\e^{-3}2^{12} & I^*_2 & 6 & 39 \\
 0 & -1+\e & 0 & -86+41\e & -297+162\e & +\e^22^{15} & I^*_5 & 6 & 38 \\
 0 & -1+\e & 0 & -262-423\e & 3063+4914\e & +\e^{-3}2^{16} & I^*_6 & 6 & 9 \\
 0 & -1+\e & 0 & -342+213\e & -2693+1666\e & -\e^{-5}2^{15} & I^*_5 & 6 & 15 \\
 0 & -1-\e & 0 & -2+\e & 1 & +2^{10} & I^*_{0} & 6 & 14 \\
 0 & -1-\e & 0 & -2-19\e & -11-24\e & +\e^{-3}2^{14} & I^*_4 & 6 & 20 \\
 0 & -1-\e & 0 & -22+\e & -11+44\e & +2^{14} & I^*_4 & 6 & 34 \\
 0 & -1-\e & 0 & -22+21\e & 57-44\e & -\e^32^{14} & I^*_4 & 6 & 36 \\
 0 & -1-\e & 0 & -86-127\e & -459-756\e & +\e^42^{15} & I^*_5 & 6 & 21 \\
 0 & -1-\e & 0 & -214+129\e & 1429-884\e & +\e^{-4}2^{15} & I^*_5 & 6 & 38 \\
 0 & \e & 0 & -2\e & -2-2\e & +\e^{-1}2^7 & II & 7 & 18 \\
 0 & \e & 0 & 1+\e & 1+2\e & -\e^62^8 & III & 7 & 5 \\
 0 & \e & 0 & 1-\e & -1 & +\e^{-1}2^8 & III & 7 & 6 \\
 0 & \e & 0 & 3+3\e & -5-10\e & -\e^62^{14} & III^* & 7 & 5 \\
 0 & \e & 0 & 3-5\e & 3-2\e & +\e^{-1}2^{14} & III^* & 7 & 6 \\
 0 & \e & 0 & -1+\e & 1 & -\e2^8 & III & 7 & 17 \\
 0 & \e & 0 & -1-9\e & 7+6\e & +\e^{-1}2^{13} & I^*_2 & 7 & 18 \\
 0 & \e & 0 & -2 & -2\e & +\e^{-2}2^7 & II & 7 & 10 \\
 0 & \e & 0 & -2-2\e & -2-4\e & +\e^62^7 & II & 7 & 3 \\
 0 & \e & 0 & -5+3\e & -5-2\e & -\e2^{14} & III^* & 7 & 17 \\
 0 & \e & 0 & -9-\e & -1+6\e & +\e^{-2}2^{13} & I^*_2 & 7 & 10 \\
 0 & \e & 0 & -9-9\e & 7+14\e & +\e^62^{13} & I^*_2 & 7 & 3 \\
 0 & \e & 0 & -11+7\e & 19-12\e & -\e^{-5}2^8 & III & 7 & 2 \\
 0 & \e & 0 & -13+8\e & 24-15\e & +\e^{-2}2^7 & II & 7 & 7 \\
 0 & \e & 0 & -45+27\e & -125+78\e & -\e^{-5}2^{14} & III^* & 7 & 2 \\
 0 & \e & 0 & -53+31\e & -161+98\e & +\e^{-2}2^{13} & I^*_2 & 7 & 7 \\
 0 & \e & 0 & -75-117\e & 419+680\e & +\e^52^8 & III & 7 & 41 \\
 0 & \e & 0 & -301-469\e & -3821-6210\e & +\e^52^{14} & III^* & 7 & 41 \\
 0 & \e & 0 & -384+237\e & -3325+2055\e & -\e^{-5}2^7 & II & 7 & 42 \\
 0 & \e & 0 & -1537+947\e & 27547-17030\e & -\e^{-5}2^{13} & I^*_2 & 7 & 42 \\
 0 & -\e & 0 & -2\e & 2+2\e & +\e^{-1}2^7 & II & 7 & 18 \\
 0 & -\e & 0 & 1+\e & -1-2\e & -\e^62^8 & III & 7 & 5 \\
 0 & -\e & 0 & 1-\e & 1 & +\e^{-1}2^8 & III & 7 & 6 \\
 0 & -\e & 0 & 3+3\e & 5+10\e & -\e^62^{14} & III^* & 7 & 5 \\
\end{array}
$$

$$
\begin{array}{cccccclrc}
 a_1 & a_2 & a_3 & a_4 & a_6 & \Delta & \mbox{type} & f & \mbox{code} \\
 0 & -\e & 0 & 3-5\e & -3+2\e & +\e^{-1}2^{14} & III^* & 7 & 6 \\
 0 & -\e & 0 & -1+\e & -1 & -\e2^8 & III & 7 & 17 \\
 0 & -\e & 0 & -1-9\e & -7-6\e & +\e^{-1}2^{13} & I^*_2 & 7 & 18 \\
 0 & -\e & 0 & -2 & 2\e & +\e^{-2}2^7 & II & 7 & 10 \\
 0 & -\e & 0 & -2-2\e & 2+4\e & +\e^62^7 & II & 7 & 3 \\
 0 & -\e & 0 & -5+3\e & 5+2\e & -\e2^{14} & III^* & 7 & 17 \\
 0 & -\e & 0 & -9-\e & 1-6\e & +\e^{-2}2^{13} & I^*_2 & 7 & 10 \\
 0 & -\e & 0 & -9-9\e & -7-14\e & +\e^62^{13} & I^*_2 & 7 & 3 \\
 0 & -\e & 0 & -11+7\e & -19+12\e & -\e^{-5}2^8 & III & 7 & 2 \\
 0 & -\e & 0 & -13+8\e & -24+15\e & +\e^{-2}2^7 & II & 7 & 7 \\
 0 & -\e & 0 & -45+27\e & 125-78\e & -\e^{-5}2^{14} & III^* & 7 & 2 \\
 0 & -\e & 0 & -53+31\e & 161-98\e & +\e^{-2}2^{13} & I^*_2 & 7 & 7 \\
 0 & -\e & 0 & -75-117\e & -419-680\e & +\e^52^8 & III & 7 & 41 \\
 0 & -\e & 0 & -301-469\e & 3821+6210\e & +\e^52^{14} & III^* & 7 & 41 \\
 0 & -\e & 0 & -384+237\e & 3325-2055\e & -\e^{-5}2^7 & II & 7 & 42 \\
 0 & -\e & 0 & -1537+947\e & -27547+17030\e & -\e^{-5}2^{13} & I^*_2 & 7 & 42 \\
 0 & 1 & 0 & 1 & 1 & -2^8 & III & 7 & 5 \\
 0 & 1 & 0 & 3 & -5 & -2^{14} & III^* & 7 & 5 \\
 0 & 1 & 0 & -1-2\e & -1-2\e & +\e^52^8 & III & 7 & 1 \\
 0 & 1 & 0 & -2 & -2 & +2^7 & II & 7 & 3 \\
 0 & 1 & 0 & -3+2\e & -3+2\e & -\e^{-5}2^8 & III & 7 & 17 \\
 0 & 1 & 0 & -5-8\e & 3+8\e & +\e^52^{14} & III^* & 7 & 1 \\
 0 & 1 & 0 & -9 & 7 & +2^{13} & I^*_2 & 7 & 3 \\
 0 & 1 & 0 & -13+8\e & 11-8\e & -\e^{-5}2^{14} & III^* & 7 & 17 \\
 0 & 1 & 0 & -33-42\e & 103+158\e & +\e^{-1}2^8 & III & 7 & 41 \\
 0 & 1 & 0 & -75+42\e & 261-158\e & -\e2^8 & III & 7 & 40 \\
 0 & 1 & 0 & -133-168\e & -957-1432\e & +\e^{-1}2^{14} & III^* & 7 & 41 \\
 0 & 1 & 0 & -301+168\e & -2389+1432\e & -\e2^{14} & III^* & 7 & 40 \\
 0 & 1+\e & 0 & -\e & -1-2\e & +\e^52^8 & III & 7 & 6 \\
 0 & 1+\e & 0 & -2+3\e & 1 & -\e^{-5}2^8 & III & 7 & 13 \\
 0 & 1+\e & 0 & -2-\e & 1 & +\e^{-1}2^8 & III & 7 & 23 \\
 0 & 1+\e & 0 & -2-2\e & -4-6\e & +\e^42^7 & II & 7 & 10 \\
 0 & 1+\e & 0 & -2-4\e & -6-10\e & +\e^52^7 & II & 7 & 18 \\
 0 & 1+\e & 0 & -2-7\e & -1 & +\e^52^{14} & III^* & 7 & 6 \\
 0 & 1+\e & 0 & -3-\e & \e & +\e^{-4}2^7 & II & 7 & 29 \\
 0 & 1+\e & 0 & -4+3\e & -5+2\e & -\e2^8 & III & 7 & 2 \\
 0 & 1+\e & 0 & -5+3\e & -6+3\e & +\e^42^7 & II & 7 & 7 \\
 0 & 1+\e & 0 & -5+4\e & 5-3\e & +\e^{-4}2^7 & II & 7 & 24 \\
 0 & 1+\e & 0 & -7+6\e & 9-5\e & -\e^{-5}2^7 & II & 7 & 33 \\
 0 & 1+\e & 0 & -10+9\e & -9+8\e & -\e^{-5}2^{14} & III^* & 7 & 13 \\
 0 & 1+\e & 0 & -10-7\e & -25-24\e & +\e^{-1}2^{14} & III^* & 7 & 23 \\
 0 & 1+\e & 0 & -10-11\e & 11+16\e & +\e^42^{13} & I^*_2 & 7 & 10 \\
 0 & 1+\e & 0 & -10-19\e & 19+32\e & +\e^52^{13} & I^*_2 & 7 & 18 \\
 0 & 1+\e & 0 & -14-7\e & -21-36\e & +\e^{-4}2^{13} & I^*_2 & 7 & 29 \\
 0 & 1+\e & 0 & -18+9\e & 31-16\e & -\e2^{14} & III^* & 7 & 2 \\
 0 & 1+\e & 0 & -22+9\e & 35-28\e & +\e^42^{13} & I^*_2 & 7 & 7 \\
 0 & 1+\e & 0 & -22+13\e & -49+28\e & +\e^{-4}2^{13} & I^*_2 & 7 & 24 \\
 0 & 1+\e & 0 & -30+21\e & -81+52\e & -\e^{-5}2^{13} & I^*_2 & 7 & 33 \\
\end{array}
$$

$$
\begin{array}{cccccclrc}
 a_1 & a_2 & a_3 & a_4 & a_6 & \Delta & \mbox{type} & f & \mbox{code} \\
 0 & 1+\e & 0 & -58-88\e & -358-574\e & +\e^{-1}2^7 & II & 7 & 43 \\
 0 & 1+\e & 0 & -147+90\e & 785-485\e & -\e2^7 & II & 7 & 42 \\
 0 & 1+\e & 0 & -234-355\e & 2275+3648\e & +\e^{-1}2^{13} & I^*_2 & 7 & 43 \\
 0 & 1+\e & 0 & -590+357\e & -6513+4004\e & -\e2^{13} & I^*_2 & 7 & 42 \\
 0 & 1-\e & 0 & \e & -1 & -\e2^8 & III & 7 & 13 \\
 0 & 1-\e & 0 & -\e & 1 & +\e^{-1}2^8 & III & 7 & 1 \\
 0 & 1-\e & 0 & -2 & -2+2\e & +\e^22^7 & II & 7 & 24 \\
 0 & 1-\e & 0 & -2+2\e & -4+2\e & -\e2^7 & II & 7 & 33 \\
 0 & 1-\e & 0 & -2+5\e & 1+2\e & -\e2^{14} & III^* & 7 & 13 \\
 0 & 1-\e & 0 & -2-3\e & -7+2\e & +\e^{-1}2^{14} & III^* & 7 & 1 \\
 0 & 1-\e & 0 & -4-7\e & 7+12\e & +\e^52^8 & III & 7 & 23 \\
 0 & 1-\e & 0 & -5-8\e & 9+15\e & +\e^22^7 & II & 7 & 29 \\
 0 & 1-\e & 0 & -10+\e & 5-6\e & +\e^22^{13} & I^*_2 & 7 & 24 \\
 0 & 1-\e & 0 & -10+9\e & 13-6\e & -\e2^{13} & I^*_2 & 7 & 33 \\
 0 & 1-\e & 0 & -18-27\e & -47-78\e & +\e^52^{14} & III^* & 7 & 23 \\
 0 & 1-\e & 0 & -22-31\e & -63-98\e & +\e^22^{13} & I^*_2 & 7 & 29 \\
 0 & 1-\e & 0 & -147-237\e & -1270-2055\e & +\e^52^7 & II & 7 & 43 \\
 0 & 1-\e & 0 & -192+117\e & 1099-680\e & -\e^{-5}2^8 & III & 7 & 40 \\
 0 & 1-\e & 0 & -590-947\e & 10517+17030\e & +\e^52^{13} & I^*_2 & 7 & 43 \\
 0 & 1-\e & 0 & -770+469\e & -10031+6210\e & -\e^{-5}2^{14} & III^* & 7 & 40 \\
 0 & -1 & 0 & 1 & -1 & -2^8 & III & 7 & 5 \\
 0 & -1 & 0 & 3 & 5 & -2^{14} & III^* & 7 & 5 \\
 0 & -1 & 0 & -1-2\e & 1+2\e & +\e^52^8 & III & 7 & 1 \\
 0 & -1 & 0 & -2 & 2 & +2^7 & II & 7 & 3 \\
 0 & -1 & 0 & -3+2\e & 3-2\e & -\e^{-5}2^8 & III & 7 & 17 \\
 0 & -1 & 0 & -5-8\e & -3-8\e & +\e^52^{14} & III^* & 7 & 1 \\
 0 & -1 & 0 & -9 & -7 & +2^{13} & I^*_2 & 7 & 3 \\
 0 & -1 & 0 & -13+8\e & -11+8\e & -\e^{-5}2^{14} & III^* & 7 & 17 \\
 0 & -1 & 0 & -33-42\e & -103-158\e & +\e^{-1}2^8 & III & 7 & 41 \\
 0 & -1 & 0 & -75+42\e & -261+158\e & -\e2^8 & III & 7 & 40 \\
 0 & -1 & 0 & -133-168\e & 957+1432\e & +\e^{-1}2^{14} & III^* & 7 & 41 \\
 0 & -1 & 0 & -301+168\e & 2389-1432\e & -\e2^{14} & III^* & 7 & 40 \\
 0 & -1+\e & 0 & \e & 1 & -\e2^8 & III & 7 & 13 \\
 0 & -1+\e & 0 & -\e & -1 & +\e^{-1}2^8 & III & 7 & 1 \\
 0 & -1+\e & 0 & -2 & 2-2\e & +\e^22^7 & II & 7 & 24 \\
 0 & -1+\e & 0 & -2+2\e & 4-2\e & -\e2^7 & II & 7 & 33 \\
 0 & -1+\e & 0 & -2+5\e & -1-2\e & -\e2^{14} & III^* & 7 & 13 \\
 0 & -1+\e & 0 & -2-3\e & 7-2\e & +\e^{-1}2^{14} & III^* & 7 & 1 \\
 0 & -1+\e & 0 & -4-7\e & -7-12\e & +\e^52^8 & III & 7 & 23 \\
 0 & -1+\e & 0 & -5-8\e & -9-15\e & +\e^22^7 & II & 7 & 29 \\
 0 & -1+\e & 0 & -10+\e & -5+6\e & +\e^22^{13} & I^*_2 & 7 & 24 \\
 0 & -1+\e & 0 & -10+9\e & -13+6\e & -\e2^{13} & I^*_2 & 7 & 33 \\
 0 & -1+\e & 0 & -18-27\e & 47+78\e & +\e^52^{14} & III^* & 7 & 23 \\
 0 & -1+\e & 0 & -22-31\e & 63+98\e & +\e^22^{13} & I^*_2 & 7 & 29 \\
 0 & -1+\e & 0 & -147-237\e & 1270+2055\e & +\e^52^7 & II & 7 & 43 \\
 0 & -1+\e & 0 & -192+117\e & -1099+680\e & -\e^{-5}2^8 & III & 7 & 40 \\
 0 & -1+\e & 0 & -590-947\e & -10517-17030\e & +\e^52^{13} & I^*_2 & 7 & 43 \\
 0 & -1+\e & 0 & -770+469\e & 10031-6210\e & -\e^{-5}2^{14} & III^* & 7 & 40 \\
\end{array}
$$

$$
\begin{array}{cccccclrc}
 a_1 & a_2 & a_3 & a_4 & a_6 & \Delta & \mbox{type} & f & \mbox{code} \\
 0 & -1-\e & 0 & -\e & 1+2\e & +\e^52^8 & III & 7 & 6 \\
 0 & -1-\e & 0 & -2+3\e & -1 & -\e^{-5}2^8 & III & 7 & 13 \\
 0 & -1-\e & 0 & -2-\e & -1 & +\e^{-1}2^8 & III & 7 & 23 \\
 0 & -1-\e & 0 & -2-2\e & 4+6\e & +\e^42^7 & II & 7 & 10 \\
 0 & -1-\e & 0 & -2-4\e & 6+10\e & +\e^52^7 & II & 7 & 18 \\
 0 & -1-\e & 0 & -2-7\e & 1 & +\e^52^{14} & III^* & 7 & 6 \\
 0 & -1-\e & 0 & -3-\e & -\e & +\e^{-4}2^7 & II & 7 & 29 \\
 0 & -1-\e & 0 & -4+3\e & 5-2\e & -\e2^8 & III & 7 & 2 \\
 0 & -1-\e & 0 & -5+3\e & 6-3\e & +\e^42^7 & II & 7 & 7 \\
 0 & -1-\e & 0 & -5+4\e & -5+3\e & +\e^{-4}2^7 & II & 7 & 24 \\
 0 & -1-\e & 0 & -7+6\e & -9+5\e & -\e^{-5}2^7 & II & 7 & 33 \\
 0 & -1-\e & 0 & -10+9\e & 9-8\e & -\e^{-5}2^{14} & III^* & 7 & 13 \\
 0 & -1-\e & 0 & -10-7\e & 25+24\e & +\e^{-1}2^{14} & III^* & 7 & 23 \\
 0 & -1-\e & 0 & -10-11\e & -11-16\e & +\e^42^{13} & I^*_2 & 7 & 10 \\
 0 & -1-\e & 0 & -10-19\e & -19-32\e & +\e^52^{13} & I^*_2 & 7 & 18 \\
 0 & -1-\e & 0 & -14-7\e & 21+36\e & +\e^{-4}2^{13} & I^*_2 & 7 & 29 \\
 0 & -1-\e & 0 & -18+9\e & -31+16\e & -\e2^{14} & III^* & 7 & 2 \\
 0 & -1-\e & 0 & -22+9\e & -35+28\e & +\e^42^{13} & I^*_2 & 7 & 7 \\
 0 & -1-\e & 0 & -22+13\e & 49-28\e & +\e^{-4}2^{13} & I^*_2 & 7 & 24 \\
 0 & -1-\e & 0 & -30+21\e & 81-52\e & -\e^{-5}2^{13} & I^*_2 & 7 & 33 \\
 0 & -1-\e & 0 & -58-88\e & 358+574\e & +\e^{-1}2^7 & II & 7 & 43 \\
 0 & -1-\e & 0 & -147+90\e & -785+485\e & -\e2^7 & II & 7 & 42 \\
 0 & -1-\e & 0 & -234-355\e & -2275-3648\e & +\e^{-1}2^{13} & I^*_2 & 7 & 43 \\
 0 & -1-\e & 0 & -590+357\e & 6513-4004\e & -\e2^{13} & I^*_2 & 7 & 42 \\
 0 & 0 & 0 & 2\e & 0 & -\e^32^9 & III & 8 & 11 \\
 0 & 0 & 0 & 8\e & 0 & -\e^32^{15} & III^* & 8 & 11 \\
 0 & 0 & 0 & -2\e & 0 & +\e^32^9 & III & 8 & 11 \\
 0 & 0 & 0 & -8\e & 0 & +\e^32^{15} & III^* & 8 & 11 \\
 0 & 0 & 0 & 2 & 0 & -2^9 & III & 8 & 11 \\
 0 & 0 & 0 & 2+2\e & 0 & -\e^62^9 & III & 8 & 11 \\
 0 & 0 & 0 & 2-2\e & 0 & +\e^{-3}2^9 & III & 8 & 11 \\
 0 & 0 & 0 & 8 & 0 & -2^{15} & III^* & 8 & 11 \\
 0 & 0 & 0 & 8+8\e & 0 & -\e^62^{15} & III^* & 8 & 11 \\
 0 & 0 & 0 & 8-8\e & 0 & +\e^{-3}2^{15} & III^* & 8 & 11 \\
 0 & 0 & 0 & -2 & 0 & +2^9 & III & 8 & 11 \\
 0 & 0 & 0 & -2+2\e & 0 & -\e^{-3}2^9 & III & 8 & 11 \\
 0 & 0 & 0 & -2-2\e & 0 & +\e^62^9 & III & 8 & 11 \\
 0 & 0 & 0 & -8 & 0 & +2^{15} & III^* & 8 & 11 \\
 0 & 0 & 0 & -8+8\e & 0 & -\e^{-3}2^{15} & III^* & 8 & 11 \\
 0 & 0 & 0 & -8-8\e & 0 & +\e^62^{15} & III^* & 8 & 11 \\
 0 & \e & 0 & -3-\e & -1 & +\e^62^9 & III & 8 & 26 \\
 0 & \e & 0 & -3-3\e & 1+2\e & +\e^62^9 & III & 8 & 35 \\
 0 & \e & 0 & -3-5\e & 3+4\e & +\e^32^9 & III & 8 & 16 \\
 0 & \e & 0 & -5-7\e & 5+8\e & +\e^62^9 & III & 8 & 12 \\
 0 & \e & 0 & -13-5\e & 3-18\e & +\e^62^{15} & III^* & 8 & 26 \\
 0 & \e & 0 & -13-13\e & -21-42\e & +\e^62^{15} & III^* & 8 & 35 \\
 0 & \e & 0 & -13-21\e & -45-66\e & +\e^32^{15} & III^* & 8 & 16 \\
 0 & \e & 0 & -21-29\e & -69-114\e & +\e^62^{15} & III^* & 8 & 12 \\
\end{array}
$$

$$
\begin{array}{cccccclrc}
 a_1 & a_2 & a_3 & a_4 & a_6 & \Delta & \mbox{type} & f & \mbox{code} \\
 0 & -\e & 0 & -3-\e & 1 & +\e^62^9 & III & 8 & 26 \\
 0 & -\e & 0 & -3-3\e & -1-2\e & +\e^62^9 & III & 8 & 35 \\
 0 & -\e & 0 & -3-5\e & -3-4\e & +\e^32^9 & III & 8 & 16 \\
 0 & -\e & 0 & -5-7\e & -5-8\e & +\e^62^9 & III & 8 & 12 \\
 0 & -\e & 0 & -13-5\e & -3+18\e & +\e^62^{15} & III^* & 8 & 26 \\
 0 & -\e & 0 & -13-13\e & 21+42\e & +\e^62^{15} & III^* & 8 & 35 \\
 0 & -\e & 0 & -13-21\e & 45+66\e & +\e^32^{15} & III^* & 8 & 16 \\
 0 & -\e & 0 & -21-29\e & 69+114\e & +\e^62^{15} & III^* & 8 & 12 \\
 0 & 1 & 0 & -1-2\e & -1+2\e & +\e^{-3}2^9 & III & 8 & 16 \\
 0 & 1 & 0 & -3 & 1 & +2^9 & III & 8 & 35 \\
 0 & 1 & 0 & -3+2\e & 1-2\e & -\e^32^9 & III & 8 & 32 \\
 0 & 1 & 0 & -3-2\e & 1+2\e & +2^9 & III & 8 & 12 \\
 0 & 1 & 0 & -5+2\e & 3-2\e & +2^9 & III & 8 & 26 \\
 0 & 1 & 0 & -5-8\e & 3-24\e & +\e^{-3}2^{15} & III^* & 8 & 16 \\
 0 & 1 & 0 & -13 & -21 & +2^{15} & III^* & 8 & 35 \\
 0 & 1 & 0 & -13+8\e & -21+24\e & -\e^32^{15} & III^* & 8 & 32 \\
 0 & 1 & 0 & -13-8\e & -21-24\e & +2^{15} & III^* & 8 & 12 \\
 0 & 1 & 0 & -21+8\e & -45+24\e & +2^{15} & III^* & 8 & 26 \\
 0 & 1-\e & 0 & -8+5\e & 7-4\e & -\e^{-3}2^9 & III & 8 & 32 \\
 0 & 1-\e & 0 & -34+21\e & -111+66\e & -\e^{-3}2^{15} & III^* & 8 & 32 \\
 0 & -1 & 0 & -1-2\e & 1-2\e & +\e^{-3}2^9 & III & 8 & 16 \\
 0 & -1 & 0 & -3 & -1 & +2^9 & III & 8 & 35 \\
 0 & -1 & 0 & -3+2\e & -1+2\e & -\e^32^9 & III & 8 & 32 \\
 0 & -1 & 0 & -3-2\e & -1-2\e & +2^9 & III & 8 & 12 \\
 0 & -1 & 0 & -5+2\e & -3+2\e & +2^9 & III & 8 & 26 \\
 0 & -1 & 0 & -5-8\e & -3+24\e & +\e^{-3}2^{15} & III^* & 8 & 16 \\
 0 & -1 & 0 & -13 & 21 & +2^{15} & III^* & 8 & 35 \\
 0 & -1 & 0 & -13+8\e & 21-24\e & -\e^32^{15} & III^* & 8 & 32 \\
 0 & -1 & 0 & -13-8\e & 21+24\e & +2^{15} & III^* & 8 & 12 \\
 0 & -1 & 0 & -21+8\e & 45-24\e & +2^{15} & III^* & 8 & 26 \\
 0 & -1+\e & 0 & -8+5\e & -7+4\e & -\e^{-3}2^9 & III & 8 & 32 \\
 0 & -1+\e & 0 & -34+21\e & 111-66\e & -\e^{-3}2^{15} & III^* & 8 & 32 \\
\end{array}
$$

}
Table 2. Elliptic curves defined over ${\mathbb{Q}}\left(\sqrt5\right)$
with good reduction away from 2.

{
\ifx\undefined\bysame
\newcommand{\bysame}{\leavevmode\hbox to3em{\hrulefill}\,}
\fi

}

\end{document}